\documentclass[reqno, a4paper]{amsart}
\usepackage{amssymb}
\usepackage{graphicx}
\usepackage{amsmath}
\usepackage[dvips]{epsfig}
\usepackage{amssymb}

\theoremstyle{definition}

\makeatletter

\begin{document}

\title[ $q$-Bernoulli numbers and $q$-Bernstein polynomials]
{A Note on  $q$-Bernoulli numbers and $q$-Bernstein polynomials}


\author[T. Kim]{Taekyun Kim}
\address{Division of General Education, Kwangwoon University, Seoul 139--701,
Korea} \email{tkkim@kw.ac.kr}

\author[C. S. Ryoo]{Cheon Seoung Ryoo}
\address{Department of Mathematics,\\
          Hannam University,  Daejeon 306-791, Korea }
\email{ryoocs@$\mbox{hnu}$.kr}

\author[H. Yi]{HeungSu Yi}
\address{Department of Mathematics,
Kwangwoon University, Seoul 139--701, Korea}
\email{hsyi@kw.ac.kr}

\thanks{This study is supported in part by  the Research Grant of Kwangwoon University in 2010}

\date{June 8, 2010}


\maketitle

\begin{abstract}
The purpose of this paper is to investigate some properties of
several $q$-Bernstein type polynomials to express the bosonic
$p$-adic $q$-integral of those polynomials on $\Bbb Z_p$.
\end{abstract}

\section{INTRODUCTION}

Let $p$ be a fixed prime number. Throughout this paper,  $\Bbb
Z_p$, $\Bbb Q_p$ and $\Bbb C_p$ will  denote the ring of $p$-adic
integers, the field of $p$-adic numbers and the field of $p$-adic
completion of the algebraic closure of $\Bbb Q_p$, respectively.
Let $\Bbb N$ be the set of natural numbers and $\Bbb Z_+ =\Bbb N
\cup \{0\}$. Let $\nu_p$ be the normalized exponential valuation
of $\Bbb C_p$ with $|p|_p=p^{-\nu_p (p)} =\frac{1}{p}$. Let $q$ be
regarded as either a complex number $q \in \Bbb C$ or a $p$-adic
number $q \in \Bbb C_p.$
 If $q\in \Bbb C$, then we always  assume $|q|<1$. If $q\in \Bbb C_p$, we usually  assume that $|1-q|_p <1$.
In this paper we use the notation of $q$-number as $$[x]_q
=[x:q]=\frac{1-q^x}{1-q}.$$
 Let  $ UD(\Bbb Z_p)$ be the set of  uniformly differentiable
 functions on $ \Bbb Z_p.$
 For $f\in UD(\Bbb Z_p)$, the bosonic
$p$-adic $q$-integral on $\Bbb Z_p$ is defined by

$$
I_{q} (f)=\int_{\Bbb Z_p } f(x) d\mu_{q} (x) = \lim_{N \rightarrow
\infty} \frac{1}{[p^N]_q} \sum_{x=0}^{p^N-1} f(x)(q)^x, \quad
(\text{see [2-5]}). \eqno(1) $$
\bigskip

In [2], the Carlitz's $q$-Bernoulli numbers are inductively
defined by
 $$
 \beta_{0,q}=1 , \text{ and }  q ( q \beta +1)^k - \beta_{k,q} =
\left \{\begin{array}{ll}
1, & \text{ if } k=1, \\
0, &  \text{ if } k>1,
\end{array}   \right.
\eqno(2)
$$
 with the usual convention of replacing
 $\beta^i$ by $\beta_{i,q}$.

The Carlitz's $q$-Bernoulli polynomials are also defined by
$$ \beta_{n,q}(x)= ( q^x \beta +[x]_q)^k = \sum_{i=0}^k \binom ki q^{ix} \beta_{i, q}[x]_q^{k-i}. \eqno(3)$$

In [2], Kim proved that the Carlitz $q$-Bernoulli numbers and
polynomials are represented by $p$-adic $q$-integral as follows:
for  $n \in \Bbb Z_+$,
$$\beta_{n, q}= \int_{\mathbb{Z}_p}[x]_q^n d\mu_q(x),  \text{ and } \beta_{n, q}(x)= \int_{\mathbb{Z}_p}[x+y]_q^n d\mu_q(y).
\eqno(4)$$

The Kim's $q$-Bernstein polynomials are defined by
$$ B_{k,n}(x,q)= \binom{n}{k}[x]_q^k
[1-x]_{q^{-1}}^{n-k},  \text{ (see [1-8])}, \eqno(5)
$$
where $n, k \in \Bbb Z_+$, and $x \in [0, 1]$.

 Let $f$ be continuous functions on $[0, 1]$.
Then the linear Kim's $q$-Bernstein  operator of order $n$ for $f$
are defined by
$$ \Bbb B_{n, q}(f \mid x)=\sum_{k=0}^n f \left(
\dfrac{k}{n}\right) B_{k,n}(x, q) , \text{ (see [5])}. $$

In this paper, we consider the $p$-adic analog of the extended
Kim's  $q$-Bernstein  polynomials on $\Bbb Z_p$ and investigate
some properties of several extended Kim's  $q$-Bernstein
polynomials to express the bosonic $p$-adic $q$-integral of those
polynomials.

\section{Extended  $q$-Bernstein Polynomials}

In this section we assume that $q \in \Bbb R$ with $0<q<1$. Let
$C[0,1]$ be the set of continuous function on $[0,1]$.

 For $f\in
C[0,1]$, we consider the extended Kim's $q$-Bernstein  operator of
order $n$ as follows:

$$
\aligned
 \Bbb B_{n, q}(f \mid x_1, x_2) &=\sum_{k=0}^n  f \left( \dfrac{k}{n}\right) \binom{n}{k}
[x_1]_q^k [1-x_2]_{q^{-1}}^{n-k} \\
&=\sum_{k=0}^n  f\left( \dfrac{k}{n}\right) B_{k,n}(x_1, x_2 \mid
q).
 \endaligned
\eqno(6)
$$
 For $n, k \in \Bbb Z_+$, and $x_1, x_2  \in [0, 1]$, the extended
  Kim's  $q$-Bernstein
polynomials  of degree $n$ are defined by
$$  B_{k,n}(x_1, x_2 \mid  q) = \binom{n}{k}
[x_1]_q^k [1-x_2]_{q^{-1}}^{n-k}.\eqno(7)$$ In the special case
$x_1=x_2=x$, then $  B_{k,n}(x_1, x_2 \mid  q)=  B_{k,n}(x, q).$

From (6) and (7) we can derive the generating function for $
B_{k,n}(x_1, x_2 \mid  q)$ as follows:

$$F_q^{(k)}( x_1, x_2 \mid t)=\dfrac{(t [x_1]_q)^k exp(t [1-x_2]_{q^{-1}})}{k!}, \eqno(8)$$
where $ k \in \Bbb Z_+$  and $x_1, x_2  \in [0, 1]$.

By (8), we get
$$ \aligned  F_q^{(k)}( x_1, x_2 \mid t) &=\sum_{n=0}^\infty  \dfrac{  [x_1]_q^k [1-x_2]_{q^{-1}}^{n} }{k! n!} t^{n+k} \\
& = \sum_{n=k}^\infty \binom nk [x_1]_q^k [1-x_2]_{q^{-1}}^{n-k} \dfrac{ t^n}{ n!}  \\
&=\sum_{n=k}^\infty   B_{k,n}(x_1, x_2 \mid  q)\dfrac{ t^n}{ n!}.
\endaligned
\eqno(9)
$$
Thus, we have
$$B_{k,n}(x_1, x_2 \mid  q) =\left \{\begin{array}{ll}
 {\binom{n}{k}}[ x_1]_q^k [1-x_2]_{q^{-1}}^{n-k},   & \mbox{ if } n \geq k \\
0, &  \mbox{ if } n < k ,
\end{array} \right.
$$
for $n,  k \in \Bbb Z_+$.

It is easy to check that
$$ B_{n-k,n}(1-x_2, 1-x_1 \mid  q^{-1})=B_{k,n}(x_1, x_2 \mid  q). \eqno(10)$$

For $ 0 \leq k \leq n,$ we have
$$ \aligned  &[1-x_2]_{q^{-1}}B_{k,n-1}(x_1, x_2 \mid  q)+ [x_1]_q  B_{k-1,n-1}(x_1, x_2 \mid  q) \\
&= [1-x_2]_{q^{-1}} \binom{n-1}{k} [x_1]_q^k
[1-x_2]_{q^{-1}}^{n-k-1}+ [x_1]_q  \binom{n-1}{k-1}
[x_1]_q^{k-1} [1-x_2]_{q^{-1}}^{n-k} \\
&=  \binom{n}{k} [x_1]_q^k [1-x_2]_{q^{-1}}^{n-k} = B_{k,n}(x_1,
x_2 \mid  q).
\endaligned
\eqno(11)
$$
Therefore,  we obtain the following theorem.

\bigskip
{ \bf Theorem 1.} For $ x_1, x_2  \in [0,1]$ and $  n, k \in
\mathbb{Z}_+,$ we have
$$ [1-x_2]_{q^{-1}}B_{k,n}(x_1, x_2 \mid  q)+ [x_1]_q  B_{k-1,n}(x_1, x_2 \mid  q)= B_{k,n+1}(x_1,
x_2 \mid  q) .
$$
 \bigskip

The partial derivative of $B_{k,n}(x_1, x_2 \mid  q)$ are also
$q$-polynomials of degree $n-1$:

 $$ \dfrac{\partial }{\partial x_1}B_{k,n}(x_1, x_2 \mid  q)= \dfrac{\log
 q}{q-1} q^{x_1} n B_{k-1,n-1}(x_1, x_2 \mid  q), $$
 and
$$ \dfrac{\partial }{\partial x_2}B_{k,n}(x_1, x_2 \mid  q)= \dfrac{\log
 q}{1-q} q^{x_2} n B_{k,n-1}(x_1, x_2 \mid  q). $$
Therefore, we obtain the following lemma.

\bigskip
{ \bf Lemma 2.} For  $  k \in \mathbb{Z}_+$  and $  n \in
\mathbb{N},   x_1, x_2  \in [0,1]$,  we have
$$ \dfrac{\partial }{\partial x_1}B_{k,n}(x_1, x_2 \mid  q)= \dfrac{\log
 q}{q-1}  n\{ (q-1) [x_1]_q  B_{k-1,n-1}(x_1, x_2 \mid  q)+  B_{k-1,n-1}(x_1, x_2 \mid  q)\}, $$
 and
$$ \dfrac{\partial }{\partial x_2}B_{k,n}(x_1, x_2 \mid  q)= \dfrac{\log
 q}{1-q}n \{ (q-1) [x_2]_q  B_{k,n-1}(x_1, x_2 \mid  q)+ B_{k,n-1}(x_1, x_2 \mid  q)\}. $$
 \bigskip

For $f=1$, by (6), we have
$$
\aligned
 \Bbb B_{n, q}(1 \mid x_1, x_2) &= \sum_{k=0}^n   B_{k,n}(x_1, x_2 \mid q)=\sum_{k=0}^n   \binom{n}{k}
[x_1]_q^k [1-x_2]_{q^{-1}}^{n-k} \\
&=(1+[x_1]_q -[x_2]_q)^n.
 \endaligned
\eqno(12)
$$
By (12), we see that
$$\dfrac{1}{(1+[x_1]_q -[x_2]_q)^n}\Bbb B_{n, q}(1 \mid x_1, x_2) =1 . $$

For $f(t)=t$, by (6), we get
$$
\aligned
 \Bbb B_{n, q}(t \mid x_1, x_2) &= \sum_{k=0}^n  \left( \dfrac{k}{n} \right)    [x_1]_q^k [1-x_2]_{q^{-1}}^{n-k}  \binom{n}{k}\\
 &= \sum_{k=1}^n   [x_1]_q^k [1-x_2]_{q^{-1}}^{n-k}  \binom{n-1}{k-1}\\
 &= [x_1]_q \sum_{k=0}^{n-1}  \binom{n-1}{k} [x_1]_q^k
 [1-x_2]_{q^{-1}}^{n-k-1},
 \endaligned
$$
where $  n \in \mathbb{N}$ and $  x_1, x_2  \in [0,1]$.

Thus, we have
  $$\dfrac{1}{(1+[x_1]_q -[x_2]_q)^{n+1}}\Bbb B_{n, q}(t \mid x_1, x_2) =[x_1]_q . $$

For $f(t)=t^2$, by (6), we have
$$
\aligned
  &\Bbb B_{n, q}(t^2 \mid x_1, x_2) \\
  &= \dfrac{n-1}{n} [x_1]_q^2 (
 1+[x_1]_q- [x_2]_{q})^{n-2} + \dfrac{[x_1]_q}{n}(
 1+[x_1]_q- [x_2]_{q})^{n-1}.
 \endaligned
 $$

 In the special case, $x_1=x_2=x$,
$$
\Bbb B_{n, q}(t^2 \mid x_1, x_2) = \dfrac{n-1}{n} [x]_q^2 +
\dfrac{[x]_q}{n}.  \eqno(13)$$

 From (13), we note
that
$$
 \lim_{ n \rightarrow  \infty} \Bbb B_{n, q}(t^2 \mid x, x) =[x]_q^2.
 $$
By (6), we see that

$$
\aligned
 \Bbb B_{n, q}(f \mid x_1, x_2) &=\sum_{k=0}^n  f \left( \dfrac{k}{n}\right) B_{k, n}(  x_1, x_2 \mid q ) \\
&=\sum_{k=0}^n  f \left( \dfrac{k}{n}\right) \binom{n}{k}
[x_1]_q^k  \sum_{j=0}^{n-k} \binom{n-k}{j}(-1)^j [x_2]_q^j   \\
&=\sum_{l=0}^n   \binom nl [x_2]_q^l  \sum_{k=0}^{l} \binom lk
(-1)^{l-k} f\left( \dfrac{k}{n}\right) \left(
\dfrac{[x_1]_q}{[x_2]_q} \right)^k.
 \endaligned
$$
From the definition of $B_{k, n}( x_1, x_2 \mid q)$, we have

$$
\aligned
 &\dfrac{n-k}{n} B_{k, n}( x_1, x_2 \mid q)+  \dfrac{k+1}{n} B_{k+1, n}( x_1, x_2 \mid q) \\
 &=\dfrac{(n-1)!}{k! (n-k-1)!}[x_1]_q^k [1-x_2]_{q^{-1}}^{n-k}
 +  \dfrac{(n-1)!}{k! (n-k-1)!}[x_1]_q^{k+1} [1-x_2]_{q^{-1}}^{n-k-1}\\
&= ([x_1]_q + [1-x_2]_{q^{-1}}) B_{k, n-1}( x_1, x_2 \mid q) \\
&=([x_1]_q + 1-[x_2]_{q}) B_{k, n-1}( x_1, x_2 \mid q),
 \endaligned
 \eqno(14)
$$
where $n \in \Bbb N $ and $k \in \Bbb Z_+, x_1, x_2 \in [0, 1].$

By the binomials theorem, we get
$$  B_{k, n}( x_1, x_2 \mid q)= \left(
\dfrac{[x_1]_q}{[x_2]_q} \right)^k  \sum_{l=k}^{n} \binom{l}{k}
\binom nl  (-1)^{l-k} [x_2]_q^l. $$

It is possible to write $[x_1]_q^k$ as a linear combination of $
B_{k, n}( x_1, x_2 \mid q)$ by using the degree evaluation
formulae and mathematical induction:
 $$\dfrac{1}{(1+[x_1]_q -[x_2]_q)^{n-1}} \sum_{k=1}^n \dfrac{ \binom k1}{\binom n1}  B_{k, n}( x_1, x_2 \mid q) =[x_1]_q . $$
By the same method, we get
 $$\dfrac{1}{(1+[x_1]_q -[x_2]_q)^{n-2}} \sum_{k=2}^n \dfrac{ \binom k2}{\binom n2}  B_{k, n}( x_1, x_2 \mid q) =[x_1]_q^2 . $$

Continuing this process, we obtain the following theorem.

\bigskip
{ \bf Theorem 3.} For  $  j \in \mathbb{Z}_+$  and $   x_1, x_2
\in [0,1]$,  we have  $$\dfrac{1}{(1+[x_1]_q -[x_2]_q)^{n-j}}
\sum_{k=j}^n \dfrac{ \binom kj}{\binom nj}  B_{k, n}( x_1, x_2
\mid q) =[x_1]_q^j . $$
 \bigskip

From Theorem 3, we have
$$\dfrac{1}{(1+[x_1]_q -[x_2]_q)^{n-j}}
\sum_{k=j}^n \dfrac{ \binom kj}{\binom nj}  B_{k, n}( x_1, x_2
\mid q) =\sum_{k=0}^j q^{\binom k2} \binom{x_1}{k}_q [k]_q! S_q(k,
j-k),$$ where $[k]_q!=[k]_q [k-1]_q \cdots [2]_q [1]_q $ and
$S_q(k, j-k)$ is the $q$-Stirling numbers of the second kind.

\section{  $q$-Bernstein Polynomials associated with the bosonic $p$-adic $q$-integral on $\Bbb Z_p$.}

In this section we assume that $q \in \Bbb C_p $ with $|
1-q|_p<1$. For  $ n \in \Bbb Z_+$, by (1), we get

$$\int_{\Bbb Z_p} [1-x+x_1]_{q^{-1}}^n d\mu_{q^{-1}}(x_1)=(-1)^n q^n \int_{\Bbb Z_p} [x+x_1]_{q}^n d\mu_{q}(x_1). \eqno(15)$$

From (4) and (15), we have
$$\beta_{n, q^{-1}}(1-x)=(-1)^n q^n \beta_{n, q}(x) \text{ for }n \in \Bbb Z_+. \eqno(16) $$

By (2), (3) and (16), we get
$$q^2 \beta_{n,q}(2)-(n+1)q^2+q=q(q \beta+1)^n=\beta_{n, q} \text{ if } n>1.$$
Thus, we have
$$\beta_{n,q}(2)=(n+1)-\dfrac{1}{q}+\dfrac{1}{q^2} \beta_{n,q}. \eqno(17)$$
\bigskip
By simple calculation, we see that
$$\int_{\Bbb Z_p} [1-x]_{q^{-1}}^n d\mu_{q}(x)=(-1)^n q^n \beta_{n, q}(-1)= \beta_{n, q^{-1}}(2). $$

From (15), (16) and (17), we can derive the following equation
(18).
$$\int_{\Bbb Z_p} [1-x]_{q^{-1}}^n d\mu_{q}(x)= q^2 \beta_{n, q^{-1}}+ (n+1)-q  \text{ if } n>1. \eqno(18)$$

Taking double bosonic $p$-adic $q$-integral on $\Bbb Z_p$, by
(18), we set
$$ \aligned
& \int_{\Bbb Z_p} \int_{\Bbb Z_p} B_{k,n}(x_1, x_2 \mid q)
d\mu_{q}(x_1) d\mu_{q}(x_2) \\
&=\binom nk \int_{\Bbb Z_p} [x_1]_{q}^n d\mu_{q}(x_1)\int_{\Bbb
Z_p} [1-x_2]_{q^{-1}}^{n-k} d\mu_{q}(x_2).
\endaligned
\eqno(19)
$$
Thus, we obtain the following theorem.

\bigskip
{ \bf Theorem 4.} For $x_1, x_2 \in [0,1]$ and $ n, k \in
\mathbb{Z}_+,$ we have
$$ \aligned
&  \int_{\Bbb Z_p} \int_{\Bbb Z_p} B_{k,n}(x_1, x_2 \mid q)
d\mu_{q}(x_1) d\mu_{q}(x_2) \\
&=\left \{\begin{array}{ll}
 {\binom{n}{k}}\beta_{k,q} (q^2 \beta_{n-k, q^{-1}}+ (n-k+1)-q ) ,   & \mbox{ if } n > k+1
 \\
0, &  \mbox{ if } n < k \\ \beta_{k, q}, & \mbox{ if } n = k \\
1,& \mbox{ if } n = k=0
\end{array} \right.
\endaligned
$$
\bigskip

From the $q$-symmetric properties(see Eq. (10)) of the
$q$-Bernstein polynomials, we have

$$ \aligned
& \int_{\Bbb Z_p} \int_{\Bbb Z_p} B_{k,n}(x_1, x_2 \mid q)
d\mu_{q}(x_1) d\mu_{q}(x_2) \\
&= \sum_{l=0}^k \binom kl (-1)^{k+l} \int_{\Bbb Z_p} \int_{\Bbb
Z_p} [1-x_1]_{q^{-1}}^{k-l} [1-x_2]_{q^{-1}}^{n-k}  d\mu_{q}(x_1)
d\mu_{q}(x_2)\\
&=\int_{\Bbb Z_p} [1-x_2]_{q^{-1}}^{n-k}   d\mu_{q}(x_2)  \{
1-k\int_{\Bbb Z_p} [1-x_1]_{q^{-1}}d\mu_{q}(x_1) \\
& \quad  \quad \quad + \sum_{l=0}^{k-2} \binom kl (-1)^{k+l}
\int_{\Bbb Z_p} [1-x_1]_{q^{-1}}^{k-l} d\mu_{q}(x_1)  \}  .
\endaligned
\eqno(20)
$$

For $ n, k \in \mathbb{Z}_+,$   by (20), we get

$$ \aligned
& \int_{\Bbb Z_p} \int_{\Bbb Z_p} B_{k,n}(x_1, x_2 \mid q)
d\mu_{q}(x_1) d\mu_{q}(x_2) \\
&=\int_{\Bbb Z_p} [1-x_2]_{q^{-1}}^{n-k}   d\mu_{q}(x_2)  \{ (
1-k-\dfrac{k}{[2]_q} )  \\
& \qquad \qquad +  \sum_{l=0}^{k-2} \binom kl (-1)^{k+l} ( q^2
\beta_{k-l, q^{-1}}+k-l+1-q) \}  .
\endaligned
\eqno(21)
$$

By (19) and (21), we obtain the following theorem.

\bigskip
{ \bf Theorem 5.} For $ n, k \in \mathbb{Z}_+,$ we have
$$\binom nk \beta_{k, q}= (
1-k-\dfrac{k}{[2]_q} ) + \sum_{l=0}^{k-2} \binom kl (-1)^{k+l} (
q^2 \beta_{k-l, q^{-1}}+k-l+1-q).
$$
\bigskip

Let  $ m, n, k \in \mathbb{Z}_+.$  Then we have

$$ \aligned
& \int_{\Bbb Z_p} \int_{\Bbb Z_p} B_{k,n}(x_1, x_2 \mid q)
B_{k,m}(x_1, x_2 \mid q)
d\mu_{q}(x_1) d\mu_{q}(x_2) \\
&=\binom nk \binom mk \int_{\Bbb Z_p} [x_1]_{q}^{2k}
d\mu_{q}(x_1)\int_{\Bbb
Z_p} [1-x_2]_{q^{-1}}^{n+m-2k} d\mu_{q}(x_2) \\
&= \binom nk \binom mk  \beta_{2k, q} \int_{\Bbb Z_p}
[1-x_2]_{q^{-1}}^{n+m-2k} d\mu_{q}(x_2) .
\endaligned
\eqno(22)
$$

By the $q$-symmetric properties of $q$-Bernstein polynomials, we
get

$$ \aligned
& \int_{\Bbb Z_p} \int_{\Bbb Z_p} B_{k,n}(x_1, x_2 \mid q)
B_{k,m}(x_1, x_2 \mid q)
d\mu_{q}(x_1) d\mu_{q}(x_2) \\
&=\sum_{l=0}^{2k} \binom {2k}{l} (-1)^{2k+l} \int_{\Bbb Z_p}
[1-x_1]_{q^{-1}}^{2k-l} d\mu_{q}(x_1) \int_{\Bbb
Z_p} [1-x_2]_{q^{-1}}^{n+m-2k} d\mu_{q}(x_2) \\
&= \int_{\Bbb Z_p} [1-x_2]_{q^{-1}}^{n+m-2k} d\mu_{q}(x_2) \{ 1-2k
\int_{\Bbb Z_p} [1-x_1]_{q^{-1}} d\mu_{q}(x_1) \\
& \quad \quad \quad  + \sum_{l=0}^{2k-2} \binom {2k}{l}
(-1)^{2k+l} \int_{\Bbb Z_p} [1-x_1]_{q^{-1}}^{2k-l} d\mu_{q}(x_1)
\} .
\endaligned
\eqno(23)
$$

By (22) and (23), we obtain the following theorem.

\bigskip
{ \bf Theorem 6.} For $ m, n, k \in \mathbb{Z}_+,$ we have
$$\binom nk \binom mk  \beta_{k, q}=
1-2k-\dfrac{2k}{[2]_q} + \sum_{l=0}^{2k-2} \binom {2k}{l}
(-1)^{2k+l} ( q^2 \beta_{2k-l, q^{-1}}+2k-l+1-q).
$$
\bigskip

Let $ n_1, n_2, \ldots , n_s,  k \in \mathbb{Z}_+,$ and $s \in
\Bbb N .$ Then

$$ \aligned
& \int_{\Bbb Z_p} \int_{\Bbb Z_p} \prod_{i=1}^s  B_{k,n_i}(x_1,
x_2 \mid q) d\mu_{q}(x_1) d\mu_{q}(x_2) \\
&= \left( \prod_{i=1}^s \binom {n_i}{k} \right) \int_{\Bbb Z_p}
[x_1]_{q}^{sk} d\mu_{q}(x_1) \int_{\Bbb
Z_p} [1-x_2]_{q^{-1}}^{n_1+ \cdots+ n_s -sk} d\mu_{q}(x_2) \\
&=  \prod_{i=1}^s \binom {n_i}{k} \beta_{sk, q} \int_{\Bbb Z_p}
[1-x_2]_{q^{-1}}^{n_1+ \cdots+ n_s -sk} d\mu_{q}(x_2) .
\endaligned
\eqno(24)
$$

By the binomial theorem, we get

$$ \aligned
& \int_{\Bbb Z_p} \int_{\Bbb Z_p} \prod_{i=1}^s  B_{k,n_i}(x_1,
x_2 \mid q) d\mu_{q}(x_1) d\mu_{q}(x_2) \\
&=  \sum_{l=0}^{sk} \binom {sk}{l}(-1)^{sk+l}  \int_{\Bbb Z_p}
[1-x_1]_{q^{-1}}^{sk-l} d\mu_{q}(x_1) \int_{\Bbb Z_p}
[1-x_2]_{q^{-1}}^{n_1+ \cdots+ n_s -sk} d\mu_{q}(x_2) .
\endaligned
\eqno(25)
$$

From (24) and (25), we note that

$$ \aligned
&   \left( \prod_{i=1}^s \binom {n_i}{k} \right) \beta_{sk, q}\\
&=  \sum_{l=0}^{sk} \binom {sk}{l}(-1)^{sk+l}  \int_{\Bbb Z_p}
[1-x_1]_{q^{-1}}^{sk-l} d\mu_{q}(x_1) \\
&= 1-sk \int_{\Bbb Z_p} [1-x_1]_{q^{-1}} d\mu_{q}(x_1)+
\sum_{l=0}^{sk-2} \binom {sk}{l}(-1)^{sk+l}  \int_{\Bbb Z_p}
[1-x_1]_{q^{-1}}^{sk-l} d\mu_{q}(x_1).
\endaligned
\eqno(26)
$$

By (26), we obtain the following theorem.

\bigskip
{ \bf Theorem 7.} Let  $s \in \Bbb N,  n_1, n_2, \ldots , n_s,  k
\in \mathbb{Z}_+.$ Then we have
$$ \left( \prod_{i=1}^s \binom {n_i}{k} \right) \beta_{sk, q} =
1-sk-\dfrac{sk}{[2]_q}  + \sum_{l=0}^{sk-2} \binom {sk}{l}
(-1)^{sk+l} ( q^2 \beta_{sk-l, q^{-1}}+sk-l+1-q).
$$

\end{document}